\documentclass[12pt,centertags]{amsart}
\usepackage{amsmath,amstext,amsthm,a4,amssymb,amscd}
\usepackage[mathscr]{eucal}
\usepackage{mathrsfs}
\usepackage{epsfig}
\numberwithin{equation}{section}

\newcommand{\field}[1]{\mathbb{#1}}

\newcommand{\bR}{\field{R}}
\newcommand{\bC}{\field{C}}
\newcommand{\bN}{\field{N}}

 \def\cC{\mathscr{C}}

\def\cO{\mathscr{O}}

\DeclareMathOperator{\End}{End}

\DeclareMathOperator{\Ker}{Ker}

\DeclareMathOperator{\Id}{Id}

\DeclareMathOperator{\td}{Td}

\DeclareMathOperator{\ch}{ch}

\newtheorem{thm}{Theorem}[section]

\theoremstyle{definition}

\theoremstyle{definition}

\newcommand{\be}{\begin{eqnarray}}
\newcommand{\ee}{\end{eqnarray}}
\newcommand{\ov}{\overline}

\newcommand{\comment}[1]{}

\begin{document}
\title{Superconnection and family Bergman kernels}
\date{\today}
\author{Xiaonan Ma}
\address{Centre de Math\'ematiques Laurent Schwartz, UMR 7640 du CNRS,
Ecole Polytechnique, 91128 Palaiseau Cedex,
France (ma@math.polytechnique.fr)}
\
\author{Weiping Zhang}
\address{Chern Institute of Mathematics \& LPMC, Nankai
University, Tianjin 300071, P.R. China.
(weiping@nankai.edu.cn)}

\begin{abstract} We establish an asymptotic expansion
for  families of Bergman kernels. The key idea is to use the
superconnection as in the local family index theorem.
\end{abstract}

\maketitle
%%%%%%%%%%%%%%%%%%%%%%%%%%%%%%%%%%%%%%%%%%%%%%%%%%%%%%%%%%%%%%%%%%%%%%%%

\begin{center}
{\Large Superconnexion et noyaux de Bergman en famille}
\end{center}
\vskip 3 mm
{\bf R\' esum\' e}. Nous annon\c cons des r\'esultats sur
 le d\'eveloppement asymptotique du noyau de Bergman en famille.

\vskip 5 mm

Let $W,S$ be smooth compact complex manifolds.
Let $\pi: W\to S$ be a holomorphic submersion with compact fiber $X$
and $\dim_\bC X=n$.

We will  add a subscript $\bR$ for the corresponding real objects.
Thus $TX$ is the holomorphic relative tangent bundle of $\pi$, and
$T_\bR X$ is the corresponding real vector bundle.
 Let $J^{T_\bR X}$ be the complex structure on $T_\bR X$.

Let $E$ be a holomorphic vector bundle on $W$.
Let $L$ be a holomorphic line bundle on $W$.
Let $h^L$, $h^E$ be Hermitian metrics on $L,E$.
Let $\nabla^L, \nabla^E$ be the  holomorphic Hermitian connections on
 $(L,h^L)$, $(E,h^E)$ with their curvatures $R^L$, $R^E$ respectively. Set
\be \label{0.1} \omega=\frac{\sqrt{-1}}{2\pi}R^L. \ee Then
$\omega$ is a smooth real $2$-form of complex type $(1,1)$ on $W$.

We suppose that $\omega$ defines a fiberwise K\"ahler
 form along the fiber $X$, i.e.
\begin{align}   \label{0.2}
g^{T_\bR X}(u,v) = \omega(u, J^{T_\bR X} v)
\end{align}
defines a (fiberwise) Riemannian metric on $T_\bR X$. We denote by
$h^{TX}$ the corresponding Hermitian metric on  $TX$.

%For a differential form $A$ on $S$, we will denote by $A^{(i)}$
%its component in $\Lambda ^i(T^*_\bR S)$.

Let $dv_X$ be the Riemannian volume form on $(X, g^{T_\bR X})$.

%then $dv_X= (\omega^n)^{(0)}/n!$.

By the Kodaira vanishing theorem, there exists $p_0\in \bN$ such
that $H^0(X, (L^p\otimes E)|_{X})$ forms a vector bundle, denoted
by $H^0(X, L^p\otimes E )$, on $S$ for $p>p_0$. From now on, we
always assume $p>p_0$.

By the Grothendieck-Riemann-Roch Theorem, for $p>p_0$, we have
\begin{align}\label{0.4}
\ch(H^0(X, L^p\otimes E))=  \int_X \td(TX) \ch(E) \ch(L^p) \quad
\mbox{in}\, \, H^\bullet (S,\bR).
\end{align}

The component in  $H^0 (S,\bR)$ of \eqref{0.4} is the
Hirzebruch-Riemann-Roch Theorem, and as $p\to +\infty$,
\begin{multline}\label{0.5}
\dim H^0(X, L^p\otimes E)
%= \left[\int_X \td(TX) \ch(E) \ch(L^p)\right]^{(0)}\\
={\rm rk} (E) \int_X \frac{c_1(L)^n}{n!} p^n \\
+ \int_X \Big(c_1(E) + \frac{{\rm rk} (E)}{2} c_1(TX)\Big)
\frac{c_1(L)^{n-1}}{(n-1)!} p^{n-1} + \cO(p^{n-2}).
\end{multline}

For $s\in S$, let $P_{p,s}$ be the orthogonal projection from
$\cC^\infty(X_s,(L^p\otimes E)_{X_s})$ onto $H^0(X_s, (L^p\otimes
E)|_{X_s})$. Let $P_{p,s}(x,x^\prime)$ $(x,x^\prime\in X_s, s\in
S)$ be the smooth kernel of $P_{p,s}$ with respect to
$dv_{X_s}(x')$.
%the Riemannian volume form $dv_X(x')$.
Then $P_{p,s}(x,x^\prime)$ is smooth on $s\in S$,
 and we denote it simply by $P_p(x,x')$, especially, $P_p(x,x)\in \End(E)_x$.

The results of \cite{T90}, \cite{Ze98}, \cite{Catlin99} tell us
that there exist $b_r\in \cC^\infty(X_s, \End(E|_{X_s}))$ such
that for any $k,l\in \bN$, there exists $C>0$ such that
\begin{align}\label{0.6}
\Big| \frac{1}{p^n}P_{p,s}(x,x)-\sum_{r=0}^k b_r(x) p^{-r}\Big|_{\cC^l(X_s)}
\leq C \, p^{-k-1}.
\end{align}
In \cite{Lu00}, \cite{Wang05},  Lu and Wang  show that the first
two coefficients $b_0,b_1$ coincide with the  corresponding terms
in the local Hirzebruch-Riemann-Roch Theorem, i.e. the leading
terms in  the Chern-Weil representative of $\td(TX) \ch(E)
\ch(L^p)$ with respect to the metrics  $h^{T X}, h^L, h^E$. We
refer to \cite{DLM06}, \cite{MM04a}, \cite{MM05b} for alternate
approaches as well as extensions to the symplectic case.

By \eqref{0.4}, in $H^2 (S,\bR)$, as $p\to \infty$, we have
\begin{multline}\label{0.7}
c_1(H^0(X, L^p\otimes E))
%= \left[\int_X \td(TX) \ch(E) \ch(L^p)\right]^{(2)}\\
={\rm rk} (E) \int_X \frac{c_1(L)^{n+1}}{(n+1)!} p^{n+1} \\
+ \int_X \Big(c_1(E) + \frac{{\rm rk} (E)}{2} c_1(TX)\Big)
\frac{c_1(L)^{n}}{n!} p^{n} + \cO(p^{n-1}).
\end{multline}

Now, in view of the Bismut local family  index theorem \cite{B86},
it is nature to ask whether the analogue of \eqref{0.6} still
holds on   higher degree levels , which will involve the curvature
of the vector bundle $H^0(X, L^p\otimes E)$ as in \eqref{0.7}.

 In this note, we  announce some results on the existence of such an expansion, and compute the
first two coefficients in the expansion.

To define a canonical connection on $H^0(X, L^p\otimes E)$, we
need to introduce a horizontal sub-bundle $T^HW$ of $TW$.

Let $T^H W$ be a sub-bundle of $TW$ such that
\begin{align} \label{0.8}
TW=T^HW \oplus TX.
\end{align}
Let $P^{TX}$ be the projection from $TW$ onto $TX$.
For $U\in TS$, let $U^H\in T^H W$ be the lift of $U$.

Clearly, \eqref{0.8} induces canonically a decomposition  $\Lambda (T^*_\bR W)
=\pi^*(\Lambda (T^*_\bR S))\widehat{\otimes}\Lambda (T^*_\bR X)$.
For a differential form $A$ on $W$, we will denote by $A^{(i)}$
its component in $\Lambda ^i(T^*_\bR S)\widehat{\otimes}\Lambda (T^*_\bR X)$.
Then $dv_X= (\omega^n)^{(0)}/n!$.

Let $T\in \Lambda ^2(T^*_\bR W)\otimes T_\bR X$ be the tensor
defined in the following way: for $U,V\in TS$, $X,Y\in TX$,
\begin{align}  \label{1.2}
\begin{split}
T(U^H,V^H)&:= -P^{TX}[U^H, V^H],\quad T(X,Y):=0,\\
 T(U^H,X)&:= \frac{1}{2} (g^{TX})^{-1}
(\mathcal{L}_{U ^H}g^{TX})X .
\end{split}
\end{align}
Let $\nabla^{L^p\otimes E}$ be the connection on $L^p\otimes E$ induced
by $\nabla^L,\nabla^E$.
For $U\in T_\bR S$, $\sigma\in \cC^\infty(S, H^0(X,L^p\otimes E))$, we define
\begin{align} \label{0.9}
\nabla^{H^0(X,L^p\otimes E)}_U \sigma
= P_p\nabla^{L^p\otimes E}_{U^H} P_p \sigma.
\end{align}
Then $\nabla^{H^0(X,L^p\otimes E)}$ is a holomorphic connection on
$H^0(X,L^p\otimes E)$ with curvature $R^{H^0(X,L^p\otimes E)}$,
but it  need not  be a Hermitian connection with respect to the
(usual) induced  $L^2$ metric $h^{H^0(X,L^p\otimes E)}$ on
$H^0(X,L^p\otimes E)$.

Let $\mathbf{k}\in T^*_\bR W$ be such that for $U\in T_\bR S$,
$X\in T_\bR X$,
\begin{align} \label{0.10}
\mathbf{k}(U^H)= \frac{1}{2} (L_{U^H}dv_X)/dv_X,\quad
\mathbf{k}(X)=0.
\end{align}
Then
\begin{align} \label{0.11}
\nabla^{\Ker(D_p)}_U= P_p(\nabla^{L^p\otimes E}_{U^H}+\mathbf{k}(U^H)) P_p,
\end{align}
is a canonical Hermitian connection on $(H^0(X,L^p\otimes
E),h^{H^0(X,L^p\otimes E)})$ with curvature $R^{\Ker(D_p)}$, but
it need not  be holomorphic.

Let $R^{H^0(X,L^p\otimes E)}(x,x')$, $R^{\Ker(D_p)}(x,x')$
$(x,x'\in X_s, s\in S)$ be the smooth kernel of the operator
$R^{H^0(X,L^p\otimes E)}$,  $R^{\Ker(D_p)}$ with respect to
$dv_X(x')$.
Then
\begin{align} \label{0.12}
R^{H^0(X,L^p\otimes E)}(x,x),\quad  R^{\Ker(D_p)}(x,x)\in
\Lambda^2(T^* S)\otimes \End(E_x).
\end{align}

If
 \begin{align}   \label{0.13}
T^HW = \{ u\in TW; \omega(u,\ov{X})=0 \, \, \mbox{for any}\,\,
X\in TX\},
\end{align}
then the triple $(\pi, g^{T_\bR X}, T^H W)$ defines a K\"ahler
fibration in the sense of \cite[Definition 1.4]{BGS88b}. In this
case, the connection $\nabla^{\Ker(D_p)}$ is the canonical
holomorphic connection on $(H^0(X,L^p\otimes
E),h^{H^0(X,L^p\otimes E)})$, and
 \begin{align}   \label{0.14}
\mathbf{k}=0,\quad \nabla^{\Ker(D_p)}= \nabla^{H^0(X,L^p\otimes E)}.
\end{align}

Let $\{w_i\}$ be an orthonormal frame of $(TX, h^{TX})$. Let
$\{e_i\}$ be an orthonormal frame of $(T_\bR X, g^{TX})$. Let
$\{g_\alpha\}$ be an frame of $TS$ and  $\{g^\alpha\}$ its dual
frame.

\begin{thm} \label{st0.1} There exist smooth sections
$b_{2,r}(x)\in \cC^\infty(W,\Lambda^2(T^* S)\otimes \End(E_x))$
which  are polynomials in $R^{TX}$, $T$,
$R^E$ {\rm(}and $R^L${\rm)}, their derivatives of order
$\leqslant 2r-1$ {\rm(}resp. $2r${\rm)} along the fiber $X$,
with
\begin{equation}\label{0.15}
b_{2,0}=- 2\pi \sqrt{-1}
\frac{(\omega^{n+1})^{(2)}}{(\omega^{n})^{(0)}} \Id_E,
\end{equation}
such that for any $k,l\in \bN$, there exists
$C_{k,\,l}>0$ such that for any $x\in W$, $p\in \bN$, $p>p_0$,
\begin{align}\label{0.16}
&\Big |\frac{1}{p^{n+1}}R^{H^0(X,L^p\otimes E)}(x,x)
- \sum_{r=0}^{k} b_{2,r}(x) p^{-r} \Big |_{\cC ^l(W)}
\leqslant C_{k,\,l}\: p^{-k-1}.
\end{align}
For $R^{\Ker(D_p)}(x,x)$, we have the similar expansion as
\eqref{0.16}, with the same leading term $b_{2,0}$ in \eqref{0.15},
and the corresponding $b_{2,r}(x)$ depends also the derivative
of $d\mathbf{k}$ of order $\leqslant 2r-1$ along the fiber $X$.

If \eqref{0.13} is verified, then
\begin{align}\label{0.17}
b_{2,1}=\left(\Big(\frac{1}{2} \langle R^{TX}w_i, \ov{w}_i\rangle
+ R^E + \frac{\sqrt{-1}}{4} g^\alpha \wedge \ov{g}^\beta
\Delta_X(\omega(g^H_\alpha, \ov{g}^H_\beta)) \Big)\omega^{n}\right)^{(2)}
/(\omega^{n})^{(0)},
\end{align}
here $\Delta_X=- \sum_i [(\nabla_{e_i})^2-
\nabla_{\nabla^{TX}_{e_i}e_i}]$ is the Bochner Laplacian along the
fiber $X$.
\end{thm}

If we take the trace on $E$ and integrate  along $X$ in
\eqref{0.16}, from \eqref{0.15} and \eqref{0.17}, we get a
refinement of \eqref{0.7} on the level of differential forms, in
the spirit of the Local Family Index Theorem.

From \eqref{0.15}, we get
\begin{equation}\label{0.18}
b_{2,0}= 2\pi g^\alpha\wedge \ov{g}^\beta
\left[- \sqrt{-1} \omega(g_\alpha^H,  \ov{g}_\beta^H)
-  \omega(g_\alpha^H, \ov{w}_j)\omega(\ov{g}_\beta^H, w_j)\right]\Id_E.
\end{equation}

By \eqref{0.16}, \eqref{0.18},
 the curvatures $R^{H^0(X,L^p\otimes E)}(x,x)$, $R^{\Ker(D_p)}(x,x)$
 provide
a natural approximation of the Monge-Amp\`ere operator on the
space of K\"{a}hler metrics. It must have relations with  the
existence problem of   geodesics on the space of K\"{a}hler
metrics (cf. \cite{Do99}, \cite{Mabuchi88}, \cite{Semmes92},
 \cite{PhSt05}).

The equation \eqref{0.18} gives also an  exact local asymptotic
behavior  of the curvature estimates    in \cite[\S 6]{Bern05}.

\comment{
Let $(X,\omega_0)$ be a compact K\"ahler manifold of dimension $n$,
we suppose that there exists a holomorphic Hermitian line bundle $(L, h^L)$
 such that $c_1(L,h^L)=\omega_0$.
Then the space of K\"ahler metrics in cohomology class $[\omega_0]$ is
\begin{align}  \label{0.20}
\mathcal{M}=\{\varphi: X\to \bR; c_1(L,e^{2\pi \varphi} h^L)
= \omega_0+ \sqrt{-1} \ov{\partial}\partial \varphi
\, \, \mbox{defines a K\"ahler form}\}/ \sim,
\end{align}
where $\varphi_1 \sim \varphi_2$ if and only if $\varphi_1=\varphi_2+ c$ for
some constant $c$.

For any complex manifold $S$ of dimension $1$
 with maps $\phi:S\to \mathcal{M}$. Let
$p_1, p_2$ be the natural projection from $W=X\times S$ onto $X,S$.
We have the holomorphic Hermitian line bundle
$(p_1^* L, e^{2\pi \phi_s} h^L)$ on $W$.
In this case, if we take $T^H W= p_2^* TS$,
$\varphi(s,x)=\phi_s(x)$,  then \eqref{0.18} reads as
\begin{multline}\label{0.21}
b_{2,0}= 2\pi g^1\wedge \ov{g}^1
\left[- \sqrt{-1} \omega(g_1,  \ov{g}_1)
- |\omega(g_1, \cdot)|^2_{h^{TX}_s}\right]\\
= 2\pi g^1\wedge \ov{g}^1\left[(\ov{\partial}^S \partial^S \varphi)
(g_1,  \ov{g}_1)- |(\ov{\partial}^X \partial^S \varphi)
(g_1,  \cdot)|^2_{h^{TX}_s}\right].
\end{multline}
In this case, the equation \eqref{0.21} gives an asymptotic
 exact local formula
of the curvature estimate given  in \cite[\S 6]{Bern05}.
}

To prove   Theorem \ref{st0.1}, we will use the superconnection
formation as in the local index theory. This is the main idea of
our work. An important feature  of  superconnection is that its
curvature is a second order differential operator along the fiber
$X$, while the superconnection itself involves derivatives along
the horizontal direction. This is also one of the points in the
local index theory. Now, by combining the formal power series
trick in \cite{MM04a}, we get in fact a general and algorithmic
way to compute the coefficients in the expansion. More details
will appear in \cite{MZ06b}.

$\ $

{\it Remark}. In this note,
 we have only formulated our results in the fiberwise positive holomorphic
line bundle case. Actually, the results hold also in the fiberwise
symplectic case, and we have  the off-diagonal expansion results
too.

$\ $

{\bf Acknowledgements.} The work of the  second author was
partially supported by MOEC and NNSFC. Part of work was done while
the first author was visiting Centre de Recerca Matem\`atica (CRM)
 in Barcelona during June and July, 2006. He would like to thank CRM for
hospitality.

\def\cprime{$'$} \def\cprime{$'$}
\providecommand{\bysame}{\leavevmode\hbox to3em{\hrulefill}\thinspace}
\providecommand{\MR}{\relax\ifhmode\unskip\space\fi MR }
\providecommand{\MRhref}[2]{%
  \href{http://www.ams.org/mathscinet-getitem?mr=#1}{#2}
}
\providecommand{\href}[2]{#2}

\end{document}